\documentclass[12pt]{article}

\usepackage[affil-it]{authblk}
      \usepackage{latexsym}
         \usepackage[reqno, namelimits, sumlimits]{amsmath}
         \usepackage{amssymb, amsfonts}
         \usepackage{amsthm}

 \newtheorem{theorem}{Theorem}[section]

 \newtheorem{lemma}[theorem]{Lemma}
 
 \newtheorem{remark}[theorem]{Remark}
 
 \newtheorem{cor}[theorem]{Corollary}

\title{Remarks on Liouville Type Theorems for Steady-State Navier-Stokes Equations}
\author{ G Seregin
  }
  
\affil{ University of Oxford, UK, and PDMI, RAS, Russia}
\date{ \today}

\begin{document}
\maketitle
\begin{abstract}
Liouville type theorems for the stationary Navier-Stokes equations are proven under certain assumptions. These assumptions are motivated by conditions that appear in Liouvile type theorems for the heat equations with a given divergence free drift. 
\end{abstract}
\setcounter{equation}{0}

\setcounter{equation}{0}
\section{Introduction}
Let us consider the classical Navier-Stokes system of partial differential equations describing the steady-state flow of a viscous incompressible fluid in the whole space $\mathbb R^3$:
\begin{equation}\label{nse}
u\cdot \nabla u-\Delta u=-\nabla p, \qquad {\rm div }\,u=0.
\end{equation}
Here, $u$ is  the three dimensional velocity field, while the pressure $p$ is a scalar function. The  classical Liouville type theorem for the stationary Navier-Stokes  equations can be formulated in the following way: prove that any bounded solution $u$ to  system
(\ref{nse})
is constant. 

A possible way to attack the above problem is to consider a similar one but for the linear equations with a given divergence free drift, i.e., 
\begin{equation} 
\label{driftse}b\cdot\nabla v -\Delta v=-\nabla q, \qquad {\rm div}\,v=0,\qquad {\rm div}\,b=0,
	\end{equation}
where $v$ is supposed to be a bounded vector field and $b$ is a divergence free drift satisfying additional conditions.

The linear setup is motivated by the  known Liouvile type theorems for  an equation
\begin{equation}
	\label{scalar}
	b\cdot\nabla f-\Delta f=0\qquad ({\rm div}\,b=0),
\end{equation} 
with respect to unknown bounded scalar function $f$ and a given function $b$. In many cases, the Laplacian in (\ref{scalar})  can be replaced with an elliptic operator, having measurable coefficients, but it is not the main issue here and we restrict ourselves to the Laplace operator with a drift term only. There are plenty of  interesting results related to Liouville type theorems for (\ref{scalar}) but we are going to mentioned only two of them that will be mimicked in our investigations of (\ref{nse}) or (\ref{driftse}). They are as follows, see \cite{NazUr} and  \cite{SSSZ2012} and references there:
\begin{theorem}
	\label{scalar case}
	Let $f$ be a bounded solutions to (\ref{scalar}). Assume that either 
\begin{equation}
	\label {BMO-1}b\in BMO^{-1}
\end{equation}	
or
\begin{equation}
	\label{Morrey}
	\sup\limits_{1<R<\infty}R^{1-\frac nq}\|b\|_{q,B(R)}<\infty
	\end{equation}
	for $ n/2<q\leq n$. Then $f$ is identically equal to a constant in $\mathbb R^n$.
\end{theorem}

How sharp Theorem \ref{scalar case} is unknown, it is a matter of further investigations.

In Theorem \ref{scalar case}, the  $n$-dimensional case is considered. In what follows, it is supposed mostly that $n=3$. 

\begin{remark} In the case $n=3$, by definition, $b\in BMO^{-1}$ if there exists a divergence free vector-valued field $\omega\in BMO$ such that $b={\rm rot}\,\omega$.
	
\end{remark}

Theorem \ref{scalar case}  can be proved with the help of Mozer's technique, for example. Although this technique is not applicable to system (\ref{driftse}), conditions (\ref{BMO-1}) and (\ref{Morrey}) may be a good perspective for the steady-state Navier-Stokes equations. As it has been shown in \cite{Ser2016} and \cite{Ser2016-2}, those conditions appear for the Stokes problem with a drift but with some additional restrictions.

\begin{theorem}\label{stokewithdrift} Assume that $v$ is a bounder solution to (\ref{driftse}) 
	and two additional conditions hold:
	\begin{equation}
\label{con2}	\sup\limits_{R>0} R^{\frac 12-\frac 3s}\|v\|_{s,B(R)}<\infty
\end{equation}
with $2\leq s\leq6$ and either
\begin{equation}
	\label{cond3}
	b\in BMO^{-1}
\end{equation}
or 
\begin{equation}
\label{con1}	
\sup\limits_{R>0} R^{1-\frac 3q}\|b\|_{L^{q,\infty}(B(R))}<\infty
\end{equation}
with $3/2<q\leq3$. Then $v\equiv0$ in $\mathbb R^3$.
\end{theorem} Here, $\|\cdot\|_{L^{q,\infty}}$ is the norm of the weak Lebesgue space $L^{q,\infty}$.

Comparing with the scalar case, the Liouville type theorem for system (\ref{driftse}) is true under  additional assumption
(\ref{con2}). The corresponding proof is quite simple and based on the Caccioppoli type inequality and scaling. Our aim is to show that in the case of non-linear system (\ref{nse}), the Liouville type theorem is valid for a similar type of conditions as in the scalar case.
\begin{theorem}
	\label{main result1} Assume that $u$ is a smooth solution to system (\ref{nse}) and $u={\rm rot}\,\omega$ in $\mathbb R^3$ for a divergence free function $\omega$.  Let, for $t>3$ and  $\alpha$ subject to  
	\begin{equation}
	\label{range}
\alpha>-\frac {t-3}{6(t-1)},	\end{equation}
the following condition hold:	
	\begin{equation}
\label{gBMO-1)}
K_\alpha(t):=\sup\limits_{R>0}R^\alpha\Big(\frac 1{|B(R)|}\int\limits_{B(R)}|\omega -[\omega]_{B(R)}|^tdx\Big)^{\frac 1t}<\infty, 
\end{equation}
where  $[\omega]_{B(R)}$ is the mean value of $\omega$ over the ball $B(R)$. 
Then $u\equiv 0$ in $	\mathbb R^3$.
\end{theorem}
\begin{remark}
\label{bmocond} Conditions (\ref{gBMO-1)}) and (\ref{range}) are satisfied if $\omega$ is a BMO-function and $\alpha=0$.
\end{remark}
\begin{theorem}
	\label{main result2} Assume that $u$ is a smooth solution to system (\ref{nse}) and
\begin{equation}
	\label{Morrey1}
	M_\beta(q):=\sup\limits_{R>0}R^\beta\Big(\frac 1{|B(R)|}	\int\limits_{B(R)}|u|^qdx\Big)^\frac 1q<\infty\end{equation}
with $3/2<q<3$ provided 
\begin{equation}
	\label{range2}
\beta>\frac {6q-3}{8q-6}. 
\end{equation} Then $u\equiv 0$ in $	\mathbb R^3$.\end{theorem}
\begin{remark}
\label{remmorrey}
For $\beta =1$, the latter condition is similar to (\ref{Morrey}). Indeed, it follows from (\ref{Morrey1}) that  $M_1(q)<\infty$ with $3/2<q< 3$ implies  $u\equiv 0$ in $\mathbb R^3$.
\end{remark}

Theorem \ref{main result2} can be improved for $s>2$.
\begin{theorem}
	\label{main result3} Assume that $u$ is a smooth solution to system (\ref{nse}) and
\begin{equation}
	\label{Morrey2}
	N_\gamma(q):=\sup\limits_{R>0}R^{\gamma-\frac 3q}\|u\|_{L^{q,\infty}}<\infty\end{equation}
with $2<q\leq3$ provided 
\begin{equation}
	\label{range3}
\gamma>\frac {4q-3}{6q-6}. 
\end{equation} Then $u\equiv 0$ in $	\mathbb R^3$.\end{theorem}

\begin{cor}\label{axialsym} Let $u$ be a bounded smooth solution to	(\ref{nse}), which is axially symmetric with respect to the axis $x_3$, and satisfy the additional assumption $|u(x)|\leq c/|x'|^\mu$ for any $x=(x',x_3)$ such that $|x'|>1$,	where 
\begin{equation}\label{decay}
	\mu>\frac 2{q_1}, \qquad q_1=\frac {15+\sqrt{33}}8.
\end{equation}
\end{cor}
This result has been already known for $\mu=1$, see \cite{KNSS2009}.

Let us mention another popular Liouville type problem  to show that any solution to system (\ref{nse}) is identically equal to zero,    provided  two conditions hold:
\begin{equation}\label{finite dissipation}
\int\limits_{\mathbb R^3}|\nabla u|^2dx<\infty.
\end{equation}
and
\begin{equation}
	\label{limit}
u(x)\to 0\quad {\rm as} \quad|x|\to \infty.
\end{equation}
 Unfortunately, whether this statement is true or not is still unknown.
 
 One of the best attempts  made  to solve the above or related problems  is presented in \cite{Galdi-book} where it is shown that  the assumption
\begin{equation}\label{9/2}
u\in L_\frac 92(\mathbb R^3)
\end{equation}
implies $u=0$. Recently,  a logarithmic improvement of condition (\ref{9/2}) has been established in \cite{ChWo2016}. It is interesting to notice that the statement of Remark \ref{axialsym} in general does not follow from (\ref{9/2}).


More Liouville type results can be found  in interesting  papers \cite{GilWein1978}, \cite{KNSS2009}, \cite{ChaeYoneda2013}, and \cite{Chae2014} and references there.

\setcounter{equation}{0}
\section{Proof of Main Result}
\subsection{Caccioppoli Type Inequality}
In the first part of the proof, we are going to  use similar arguments as in \cite{Ser2016} with some changes.
We fix $R>0$ arbitrarily and  take a non-negative cut-off function $\varphi\in C^\infty_0(B(R))$ with the following properties:  $\varphi(x)=1$ in $B(\varrho)$, $\varphi(x)=0$ out of $B(r)$, and $|\nabla^k \varphi(x)|\leq c/(r-\varrho)^k$ for any $R/2\leq\varrho<r\leq R$ and $k=1,2,3,4$. 

Assume also that we are given two divergence free functions $\overline \omega$ and $\widetilde \omega$ such that  $u={\rm rot}\,\overline \omega={\rm rot}\,\widetilde \omega$ in $B(R)$. 

For a given  $2<s<\infty$, there exists a  constant $c_0=c_0(s)>0$ and a function $w\in W^1_s(B(r))$, vanishing on $\partial B(r)$,  such that ${\rm div}\,w=\nabla \varphi\cdot u$ and 
\begin{equation}\label{Bogovskii}
\int\limits_{B(r)}|\nabla w|^tdx\leq c_0(s)\int\limits_{B(r)}|\nabla \varphi\cdot  u|^tdx\end{equation}
for $t=2$ and for $t=s$.

Now, we let us test the Navier-Stokes equations (\ref{nse})  with the function $\varphi  u-w$. After integration by parts over $B(r)$, we find the following identity
$$\int\limits_{B(r)}\varphi |\nabla u|^2dx=
-\int\limits_{B(r)}\nabla u :(\nabla \varphi\otimes u) dx+\int\limits_{B(r)}\nabla w :\nabla u dx+$$
$$-\int\limits_{B(r)}(u\cdot\nabla u)\cdot\varphi u dx+\int\limits_{B(r)}(u\cdot\nabla u)\cdot wdx=I_1+I_2+I_3+I_4.$$

The first two terms  can be estimated easily. As a result, 
$$|I_1|+|I_2|\leq c\Big(\int\limits_{B(r)}|\nabla u|^2dx\Big)^\frac 12\Big(\int\limits_{B(r)}|\nabla \varphi|^2| u|^2dx\Big)^\frac 12. $$

To estimate $I_3$ and $I_4$, we are going to use  integration by parts and elementary properties of the differential operator ${\rm rot}$. Indeed, we have
$$|I_3|=\Big| \int\limits_{B(r)} u_ju_{i,j} u_i\varphi dx\Big|=\frac12\Big| \int\limits_{B(r)} {\rm rot}\,\overline\omega\cdot \nabla|u|^2 \varphi dx\Big|= \frac12\Big| \int\limits_{B(r)} \overline \omega\cdot (\nabla|u|^2\times \nabla \varphi) dx\Big|$$
$$\leq \Big(\int\limits_{B(r)}|\nabla u|^2dx\Big)^\frac 12
\Big(\int\limits_{B(r)}|\nabla\varphi|^2|\overline \omega|^2| u|^2dx\Big)^\frac 12\leq
$$$$\leq\Big(\int\limits_{B(r)}|\nabla u|^2dx\Big)^\frac 12\Big(\int\limits_{B(r)}|\nabla \varphi|^s|u|^sdx\Big)^\frac 1s \Big(\int\limits_{B(r)}|\overline \omega|^{\frac{2s}{s-2}}dx\Big)^{\frac{s-2}{2s}}.$$
It remains to evaluated $I_4$:
$$|I_4|=\Big| \int\limits_{B(r)} ({\rm rot}\,\overline\omega)_j u_{k,j} w_kdx\Big|=\Big| \int\limits_{B(r)} \overline \omega\cdot(\nabla u_k\times \nabla w_k)  dx\Big|\leq $$
$$\leq \Big(\int\limits_{B(r)}|\nabla u|^2dx\Big)^\frac 12
\Big(\int\limits_{B(r)}|\overline \omega|^2|\nabla w|^2dx\Big)^\frac 12\leq $$
$$\leq \Big(\int\limits_{B(r)}|\nabla u|^2dx\Big)^\frac 12
\Big(\int\limits_{B(r)}|\nabla w|^sdx\Big)^\frac 1s\Big(\int\limits_{B(r)}|\overline \omega|^\frac {2s}{s-2}dx\Big)^\frac {s-2}{2s}\leq $$
$$\leq c(s)\Big(\int\limits_{B(r)}|\nabla u|^2dx\Big)^\frac 12
\Big(\int\limits_{B(r)}|\nabla\varphi|^s|u|^sdx\Big)^\frac 1s\Big(\int\limits_{B(r)}|\overline \omega|^\frac {2s}{s-2}dx\Big)^\frac {s-2}{2s}.
$$
So, after application of Young inequality, we arrive at  the following estimate
$$\int\limits_{B(\varrho )}|\nabla u|^2dx\leq \frac 18\int\limits_{B(r)}|\nabla u|^2dx+c\int\limits_{B(r)}|\nabla \varphi|^2| u|^2dx+$$$$
+c(s)\Big(\int\limits_{B(r)}|\nabla\varphi|^s|u|^sdx\Big)^\frac 2s\Big(\int\limits_{B(r)}|\overline \omega|^\frac {2s}{s-2}dx\Big)^\frac {s-2}{s}.$$
Our next step is the application of a multiplicative inequality
$$
\Big(\int\limits_{B(r)}|\nabla\varphi|^s|u|^sdx\Big)^\frac 2s\leq c\Big(\int\limits_{B(r)}|\nabla \varphi|^2| u|^2dx\Big)^{1-\lambda}\Big(\Big(\int\limits_{B(r)}|\nabla \varphi|^2|\nabla u|^2dx\Big)^\lambda+$$$$+\Big(\int\limits_{B(r)}|\nabla |\nabla \varphi||^2|u|^2dx\Big)^\lambda\Big)$$
with $\lambda=3\frac {s-2}{2s}$. It is legal under the additional assumption $s<6$. Applying the Young inequality one more time, we find 
$$\int\limits_{B(\varrho )}|\nabla u|^2dx\leq \frac 14\int\limits_{B(r)}|\nabla u|^2dx+c\int\limits_{B(r)}|\nabla \varphi|^2| u|^2dx+$$$$
+c(s)\Big(\int\limits_{B(r)}|\nabla \varphi|^2| u|^2dx\Big)^{1-\lambda}\Big(\int\limits_{B(r)}|\nabla |\nabla \varphi||^2|u|^2dx\Big)^\lambda\Big(\int\limits_{B(r)}|\overline \omega|^\frac {2s}{s-2}dx\Big)^\frac {s-2}{s}+$$$$+c(s)\Big(\frac 1{(r-\varphi)^2}\Big)^\frac \lambda{1-\lambda}\int\limits_{B(r)}|\nabla \varphi|^2| u|^2dx\Big(\int\limits_{B(r)}|\overline \omega|^\frac {2s}{s-2}dx\Big)^\frac {s-2}{s(1-\lambda)}.$$
To proceed further, we need the following auxiliary lemma.
\begin{lemma}
	\label{lemma} Let $\psi$ be a bounded non-negative twice differentiable functions with compact  support in $B(R)$. Then
	$$\int\psi|u|^2dx\leq c\Big(\int \psi |\nabla u|^2dx\Big)^\frac 12(\int\psi |\widetilde \omega|^2dx\Big)^\frac 12+c\int|\nabla^2 \psi|^2|\widetilde \omega|^2dx. $$
\end{lemma}
Proof of Lemma \ref{lemma} is based on the elementary identity:
$$
\int \psi|{\rm rot}\,\widetilde\omega|^2dx=\int\psi {\rm rot}\,{\rm rot}\, \widetilde\omega\cdot \widetilde\omega dx-\int\nabla^2\psi: \widetilde\omega\otimes \widetilde\omega dx+\frac 12 \int \Delta \psi|\widetilde\omega|^2dx.
$$

It follows from Lemma \ref{lemma} that
$$\int\limits_{B(r)}|\nabla \varphi|^2|u|^2dx\leq$$$$\leq \frac c{(r-\varrho)^2}\Big[\Big(\int\limits_{B(r)\setminus B(\varrho)}|\widetilde\omega|^2dx\Big)^\frac 12\Big(\int\limits_{B(r)\setminus B(\varrho)}|\nabla u|^2dx\Big)^\frac 12+ \frac 1{(r-\varrho)^2}\int\limits_{B(r)\setminus B(\varrho)}|\widetilde\omega|^2dx\Big]$$
and 
$$\int\limits_{B(r)}|\nabla^2 \varphi|^2|u|^2dx\leq$$$$ \leq \frac c{(r-\varrho)^4}\Big[\Big(\int\limits_{B(r)\setminus B(\varrho)}|\widetilde\omega|^2dx\Big)^\frac 12\Big(\int\limits_{B(r)\setminus B(\varrho)}|\nabla u|^2dx\Big)^\frac 12+ \frac 1{(r-\varrho)^2}\int\limits_{B(r)\setminus B(\varrho)}|\widetilde\omega|^2dx\Big].$$

So, the main inequality can be transformed to the following form
$$\int\limits_{B(\varrho )}|\nabla u|^2dx\leq \frac 14\int\limits_{B(r)}|\nabla u|^2dx+$$
$$+\frac c{(r-\varrho)^2}\Big[\Big(\int\limits_{B(r)\setminus B(\varrho)}|\widetilde\omega|^2dx\Big)^\frac 12\Big(\int\limits_{B(r)\setminus B(\varrho)}|\nabla u|^2dx\Big)^\frac 12+ \frac 1{(r-\varrho)^2}\int\limits_{B(r)\setminus B(\varrho)}|\widetilde\omega|^2dx\Big]\times$$
$$\times \Big(1+\Big(\frac 1{(r-\varrho)^3}\int\limits_{B(r)}|\overline\omega|^\frac{2s}{s-2}dx\Big)^{2\frac {s-2}{6-s}}+\Big(\frac 1{(r-\varrho)^3}\int\limits_{B(r)}|\overline\omega|^\frac{2s}{s-2}dx\Big)^{\frac {s-2}{s}}\Big).$$

Now, we are going to apply the Young inequality and find 
$$\int\limits_{B(\varrho )}|\nabla u|^2dx\leq \frac 12\int\limits_{B(r)}|\nabla u|^2dx+$$
$$+\frac c{(r-\varrho)^4}\int\limits_{B(r)\setminus B(\varrho)}|\widetilde\omega|^2dx\times$$
$$\times \Big(1+\Big(\frac 1{(r-\varrho)^3}\int\limits_{B(r)}|\overline\omega|^\frac{2s}{s-2}dx\Big)^{2\frac {s-2}{6-s}}+\Big(\frac 1{(r-\varrho)^3}\int\limits_{B(r)}|\overline\omega|^\frac{2s}{s-2}dx\Big)^{\frac {s-2}{s}}\Big)^2\leq $$
$$\leq
\frac 12\int\limits_{B(r)}|\nabla u|^2dx+
$$
$$+\frac c{(r-\varrho)^4}\int\limits_{B(r)\setminus B(\varrho)}|\widetilde\omega|^2 \Big(1+\Big(\frac 1{(r-\varrho)^3}\int\limits_{B(r)}|\overline\omega|^\frac{2s}{s-2}dx\Big)^{4\frac {s-2}{6-s}}\Big).$$

Using known iterative arguments, see \cite{Giaquinta1983}, we can deduce from the latter inequality the following:
\begin{equation}
	\label{Caccio}
\int\limits_{B(R/2 )}|\nabla u|^2dx\leq \frac cR\Big(\frac 1{|B(R)|}	\int\limits_{B(R)\setminus B(R/2)}|\widetilde\omega|^2dx\Big)
\Big(1+$$$$+\Big(\frac 1{|B(R)|}\int\limits_{B(R)}|\overline\omega|^\frac{2s}{s-2}dx\Big)^{4\frac {s-2}{6-s}}\Big).$$
Now, after the change 
$$t=\frac {2s}{s-2},$$
we find 
$$\int\limits_{B(R/2 )}|\nabla u|^2dx\leq \frac cR\Big(\frac 1{|B(R)|}	\int\limits_{B(R)\setminus B(R/2)}|\widetilde\omega|^2dx\Big)\times$$$$\times
\Big(1+\Big(\frac 1{|B(R)|}\int\limits_{B(R)}|\overline\omega|^tdx\Big)^{\frac {4}{t-3}}\Big),
\end{equation}
with the range for $t>3$. 

Obviously, Caccioppoli's type inequality (\ref{Caccio}) is going to be the main  sought 
 for proving Liouville type theorems. Indeed, we need to find  reasonable conditions under which the right hand side of (\ref{Caccio}) tends to zero as $R\to \infty$. 
 
 \subsection{Proof of Theorem \ref{main result1} }

Here, we let 
$$\overline\omega =\widetilde \omega=\omega- [\omega ]_{B(R)}$$
and  assume that 
and conditions (\ref{gBMO-1)}) and (\ref{range}) hold. Then 
(\ref{Caccio}) implies
\begin{equation}\nonumber	\label{important}
\int\limits_{B(R/2 )}|\nabla u|^2dx\leq \frac cR\Big(\frac 1{|B(R)|}	\int\limits_{B(R)}|\overline\omega|^2dx\Big)\times$$$$\times
\Big(1+\Big(\frac 1{|B(R)|}\int\limits_{B(R)}|\overline\omega|^tdx\Big)^{\frac {4}{t-3}}\Big),
\end{equation}
with the right hand side tending to zero as $R\to\infty$.

\subsection{Proof of Theorem \ref{main result2}} 

Here, we let $$\overline\omega =\widetilde \omega=\omega^R,
$$ where $\omega^R$ is a unique solution the following boundary value problem 
$${\rm rot}\,\omega^R=u,\quad {\rm div}\,\omega^R=0$$
in $B(R)$ and 
$$\omega^R\cdot\nu=0$$
on $\partial B(R)$, where $\nu $ is a unit outward normal to $\partial B(R)$. Then (\ref{Caccio}) implies
\begin{equation}	\label{important1}
\int\limits_{B(R/2 )}|\nabla u|^2dx\leq \frac cR\Big(\frac 1{|B(R)|}	\int\limits_{B(R)}|\omega^R|^2dx\Big)\times$$$$\times
\Big(1+\Big(\frac 1{|B(R)|}\int\limits_{B(R)}|\omega^R|^tdx\Big)^{\frac {4}{t-3}}\Big),
\end{equation}

By scaling, we have the following inequalities, see \cite{BB1974},
$$\int\limits_{B(R)}|\omega^R|^2dx
\leq cR^2\Big(\int\limits_{B(R)}|{\rm rot}\,\omega^R|^qdx\Big)^\frac 2qR^{3(1-\frac 2q)}$$
and
$$\Big(\int\limits_{B(R)}|\omega^R|^{\frac {3q}{3-q}}dx\Big)^\frac {3-q}{3q}\leq \Big(\int\limits_{B(R)}|{\rm rot}\,\omega^R|^qdx\Big)^\frac 1q$$
Next, letting $t=3q/(3-q)$, we can re-write (\ref{important1})  so that
$$\int\limits_{B(R/2 )}|\nabla u|^2dx\leq  cR\Big(\frac 1{|B(R)|}	\int\limits_{B(R)}|u|^qdx\Big)^\frac 2q\times$$$$\times
\Big(1+R^\frac {4q}{2q-3}\Big(\frac 1{|B(R)|}\int\limits_{B(R)}|u|^qqx\Big)^{\frac {4}{2q-3}}\Big).
$$
Letting $t=3q/(3-q)$, we have 
$$\int\limits_{B(R/2 )}|\nabla u|^2dx\leq cM_\beta(q)\Big[R^{1-2\beta}+(M_\beta(q))^\frac{4q}{2q-3}R^{\frac {4q}{2q-3}(1-\beta)+(1-2\beta)}\Big]\to 0
$$
as $R\to\infty$ provided (\ref{range2}) holds.

 \subsection{Proof of Theorem \ref{main result3} }
 Repeating arguments from (\cite{Ser2016-2}), we have the following estimate
 $$\int\limits_{B(R/2)}|\nabla u|^2dx\leq 
 \frac c{R^2}\int\limits_{B(R)}|u-[u]_{B(R)}|^2dx\Big[1+\Big(R^{1-\frac 3q}\|u\|_{L^{q,\infty}}\Big)^{\frac{2q}{2q-3}}\Big]. $$
After application of H\"older inequality, we find 
$$\int\limits_{B(R/2)}|\nabla u|^2dx\leq \frac cR(N_\gamma(q)R^{1-\gamma})^2[1+(N_\gamma(q)R^{1-\gamma})^{\frac {2q}{2q-3}}].$$
The right hand side of the latter inequality tends to  zero as $R\to\infty$ provided (\ref{range3}) holds.

 \subsection{Proof of Corollary \ref{axialsym}}

Here, we are going to replace balls $B(R)$ with cylinders $C(R)=
b(R)\times ]-R,R[$, where $b(R)$ is the two-dimensional ball of radius $R$ centred at the origin.
We also can replace the weak Lebesgue space $L^{q,\infty}(B(R))$ with the usual Lebesgue space $L_q(B(R))$. Without loss of generality, we may assume that $|u|\leq 1$.

Assume that condition (\ref{decay}) holds. We can find $q<q_1$ so that $\mu>2/q>2/q_1$. This allows us to set $\gamma=\frac 2q.$ By the choice of $q$, we have 
$$\gamma>\frac {4q-3}{6q-6}$$
and thus restriction   (\ref{range3}) is fulfilled.

  Next, it is obvious 
that 
$N_\gamma(q)=max\{I_1,I_2\}$, where 
$$I_1=\sup\limits_{0<R\leq 1}R^{-\frac 1q}\Big(\int\limits_{C(R)}|u|^qdx\Big)^\frac 1q\leq c<\infty$$
and 
$$\sup\limits_{R> 1}R^{-\frac 1q}\Big(\int\limits_{C(R)}|u|^qdx\Big)^\frac 1q\leq 
J_1+J_2.$$
Here,
$$J_1=\sup\limits_{R> 1}R^{-\frac 1q}\Big(\int\limits_{b(1)\times ]-R,R[}|u|^qdx\Big)^\frac 1q, $$$$ J_2=\sup\limits_{R> 1}R^{-\frac 1q}\Big(\int\limits_{(b(R)\setminus b(1))\times ]-R,R[}|u|^qdx\Big)^\frac 1q.$$
Introducing polar coordinates, we observe that $J_1$ is bounded and 
$$J_2\leq c\sup\limits_{R> 1}R^{-\frac 1q}\Big(2R\int\limits^R_1\varrho^{1-\mu q}d\varrho 
\Big)^\frac 1q\leq c\sup\limits_{R> 1}R^{-\frac 1q}\Big(2R(1-R^{2-\mu q})\Big)^\frac 1q\leq c.$$


\begin{thebibliography}{99}

\bibitem{BB1974}
Bourguignon, J.P., Brezis, H., Remarks on Euler Equation, J. Funct. Anal., 15, 341-363 (1974).


\bibitem{Chae2014}
Chae, D., Liouville-Type Theorem for the Forced Euler Equations and the Navier-Stokes Equations.
Commun. Math. Phys.326: 37-48 (2014).
\bibitem{ChaeYoneda2013}
Chae, D.,  Yoneda, T., On the Liouville theorem for the stationary Navier-Stokes equations in a
critical space, J. Math. Anal. Appl. 405 (2013), no. 2, 706-710.

\bibitem{ChWo2016}
Chae, G., Wolf, J.,
On Liouville type theorems for the steady Navier-Stokes equations in $R^3$, arXiv:1604.07643.

\bibitem{Galdi-book}
 Galdi, G. P. An introduction to the mathematical theory of the Navier-Stokes equations. Steady-state problems. Second edition. Springer Monographs in Mathematics. Springer, New York, 2011. xiv+1018 pp.
 
 \bibitem{Giaquinta1983}
 Giaquinta, M., Multiple integrals in the calculus of variations and nonlinear elliptic systems. Annals of Mathematics Studies, 105. Princeton University Press, Princeton, NJ, 1983. vii+297 pp. 

 \bibitem{GilWein1978}
  Gilbarg, D.,  Weinberger, H. F. Asymptotic properties of steady plane solutions of the Navier-
Stokes equations with bounded Dirichlet integral, Ann. Scuola Norm. Sup. Pisa Cl. Sci.(4) 5
(1978), no. 2, 381-404. 

\bibitem{KNSS2009}
Koch, G., Nadirashvili, N., Seregin, G., Sverak, V.,
Liouville theorems for the Navier-Stokes equations and applications,
Acta Mathematica, 203 (2009), 83--105.

\bibitem{NazUr}
Nazarov, A.~I., Uraltseva, N.~N.,
The Harnack inequality and related properties for solutions to elliptic and parabolic equations with divergence-free lower-order coefficients, 
St. Petersburg Mathematical Journal, 2012, 23:1, 93--115.





\bibitem{Ser2016}
Seregin, G., Liouville type theorem for stationary
Navier-Stokes equations, Nonlinearity, 29 (2016), 2191--2195.
\bibitem{Ser2016-2}
Seregin, G., A Liouville Type Theorem for Steady-State Navier-Stokes Equations, arXiv:1611.01563

\bibitem{SSSZ2012}
Seregin, G., Silvestre, L., Sverak, V., Zlatos, A.,
 On divergence-free drifts,  Journal of Differential Equations, Volume 252, Issue 1, January 2012, pp.  505-540.








 


\bibitem{Stein1970}
 Stein, Elias M. Singular integrals and differentiability properties of functions. Princeton Mathematical Series, No. 30 Princeton University Press, Princeton, N.J. 1970 xiv+290 pp. 
 
\end{thebibliography}
\end{document}